\documentclass[11pt,english,reqno]{amsart}
\usepackage{times}
\usepackage[T1]{fontenc}
\usepackage{babel}
\usepackage{amssymb}

\makeatletter

\providecommand{\LyX}{L\kern-.1667em\lower.25em\hbox{Y}\kern-.125emX\@}

\theoremstyle{plain}
\newtheorem*{thm*}{Theorem}
\theoremstyle{plain}
\newtheorem{lem}{Lemma}
\theoremstyle{plain}
\newtheorem*{prop*}{Proposition}
\theoremstyle{definition}
\newtheorem*{defn*}{Definition}

\newcommand{\ZZ}{\mathbb {Z}}

\newcommand{\RR}{\mathbb {R}}

\newcommand{\spec}{\mathrm{spec}\,}

\newcommand{\half}{{\textstyle \frac{1}{2}}}

\newcommand{\brk}{\\ }

\def\th@nopoint{
  \thm@headpunct{} 
  \itshape 
}

\renewcommand{\@makefnmark}{\mbox{$^{}$}}

\makeatother
\begin{document}

\title{Random Menshov Spectra}

\author{Gady Kozma and Alexander Olevski\v{\i}}



\begin{abstract}
We show that a spectrum \( \Lambda \) of
frequencies \( \lambda \) obtained
by  a random perturbation of the integers allows one to represent any
measurable function \( f \) on \( \RR  \) by an almost everywhere
converging sum of harmonics: \[
f=\sum _{\Lambda}c_{\lambda}e^{i\lambda t}\]
almost surely.
\end{abstract}
\maketitle

\section{Introduction}

This paper concerns the representation of functions
by  series of 
\footnote{\rm Research
supported in part by the Israel Science Foundation.

2000 Mathematics Subject Classification: 42A63,
42A61, 42A55.

Keywords and phrases: random spectra;
representation of functions by trigonometric
series.}
exponentials
which converge almost everywhere (a.e.). According to Menshov's theorem
(1941, see \cite{bary}) every \( 2\pi  \)-periodic measurable function
\( f \) admits a representation as\begin{equation}
\label{olv1}
f(t)=\sum _{k\in \ZZ }c(k)e^{ikt}\quad
\textrm{a}.\textrm{e}.
\end{equation}
Among the many generalizations and analogs of this fundamental result,
there exists a version for the non-periodic case: Davtjan \cite{Davtjan}
proved that the corresponding representation on \( \RR  \) can be
obtained if instead of the sum over integers one considers a 
{}``trigonometric
integral'' which involves all real frequencies.

In our  paper \cite{KO-99} it was shown that most of the frequencies
are redundant. Namely, by appropriate small perturbations of the integers
we constructed a spectrum of frequencies \(
H=\{\lambda(k),\, k\in \ZZ \} \)
such that any \( f\in L^{0}(\RR ) \) (that is, a measurable function
on \( \RR  \)) can be decomposed as\begin{equation}
\label{olv2}
f(t)=\sum _{k\in \ZZ }c(k)e^{i\lambda(k)t}\quad \textrm{a}.\textrm{e}.
\end{equation}

The aim of this note is to show that this is  not an exceptional feature
of the constructed spectrum. In fact, by choosing the perturbations
randomly one gets such a property with probability 1.

Analogously to the periodic case (see \cite{KO01}) we introduce
the following

\begin{defn*}
A sequence \( \Lambda=\{\lambda(k),\,
...<\lambda(-1)<\lambda(0)<\lambda(1)<...\} \) is called
a Menshov spectrum for \( \RR  \) if for any \( f\in L^{0}(\RR ) \)
there are coefficients \( \{c(k)\} \) such that the decomposition
(\ref{olv2}) holds (convergence in (\ref{olv2}) is
 understood in the sense of
the  limit of symmetric partial sums, i.e. \( \lim _{x\rightarrow \infty
}\sum 
_{|\lambda|<x} \)).
\end{defn*}
\begin{thm*}
Let \( r(n) \) be independent variables uniformly distributed on
the segment \( [-\half ,\half ] \). Then the sequence \begin{equation}
\label{olv3}
\lambda(n)=n+r(n),\quad n\in \ZZ
\end{equation}
is almost surely a Menshov spectrum for \( \RR  \).
\end{thm*}
The proof is based on the technique used in our recent paper \cite{KO01}.
We refer the reader to this paper for some historical comments and
additional references.

\section{Preliminaries}

By a trigonometric polynomial \( P \) we mean a finite linear combination
of exponentials with real (not necessarily integer) frequencies \( 
...<\lambda(-1)<\lambda(0)<... \).
We call the set of \( \lambda \)'s involved, the spectrum of \( P \) and
denote it by \( \spec P \). The corresponding coefficients are denoted
as \( \widehat{P}(\lambda) \), so\[
P=\sum _{\spec P}\widehat{P}(\lambda)e^{i\lambda
x},\quad x\in \RR\ .\]

We denote\[
\deg P=\max _{\spec P}|\lambda|\ .\]
As usual\[
||\widehat{P}||_{1}:=\sum |\widehat{P}(\lambda)|,\qquad ||\widehat{P}||_{\infty 
}:=\max |\widehat{P}(\lambda)|\ .\]
Let \( P^{*} \) be the (non-symmetric) majorant of partial sums:\[
P^{*}(x):=\max _{a<b}\left| \sum _{\spec P\cap 
[a,b]}\widehat{P}(\lambda)e^{i\lambda x}\right| \ .\]
For a given \( P \) and \( l\in \ZZ ^{+} \) we denote by \( P_{[l]} \)
the {}``contracted'' polynomial:\[
P_{[l]}(x)=P(lx)\ .\]
The following {}``special products'' are used:\[
H=Q_{[l]}P\ .\]
If \( \spec Q\subset \ZZ  \) and \( l>2\deg P \) then this product
has a simple structure which provides the following
estimate (compare
with (10) in \cite{KO01}):\begin{equation}
\label{olv4}
H^{*}(x)\leq \left| P(x)\right| \cdot ||Q^{*}||_{L^{\infty }(-\pi ,\pi 
)}+2P^{*}(x)\cdot ||\widehat{Q}||_{\infty }\ .
\end{equation}
We will use the following

\begin{lem}
\label{Lemma_Korner}(see \cite{KO01}, Lemma 2.1) Given any \( \epsilon >0 
\),
\( \delta >0 \), there exists a trigonometric polynomial \( P=P_{\epsilon 
,\delta } \)
with integer spectrum such that
\end{lem}
\begin{enumerate}
\item \label{req_Qhat_lem_korner}\( \widehat{P}(0)=0 \), \( 
||\widehat{P}||_{\infty }<\delta  \);
\item \label{req_Q1_lem_korner}\( \mathbf{m}\left\{ x\in [-\pi ,\pi ]\, :\, 
|P(x)-1|>\delta \right\} <\epsilon  \);
\item \label{req_Sq_lem_korner}\( ||P^{*}||_{\infty }<C\epsilon ^{-1} \).
\end{enumerate}

\section{Proof of the theorem}

\subsection{{}}

The result is an easy consequence of the following (nonstochastic)

\begin{prop*}
Let \( \Lambda=\{\lambda(n)\} \), \( \lambda(n)=n+r(n) \), \( n\in \ZZ  \), and suppose
that for every \( k\in \ZZ ^{+} \) there exists a number \( l=l(k) \)
in \( \ZZ ^{+} \) such that\begin{equation}
\label{olv5}
|r(sl+q)-2^{-|q|+1}|<\frac{1}{k^{2}},\quad
0<|s|<k,\, |q|<k\ .
\end{equation}
Then \( \Lambda \) is a Menshov spectrum for \( \RR  \).
\end{prop*}
To deduce the theorem from the proposition it is enough to mention
that if we fix \( k \) and run \( l \) over a sufficiently fast
increasing sequence \( \{l_{j}\} \) then the events \( B_{j} \)
that the inequalities above are fulfilled for \( l=l_{j} \) are mutually
independent and each has a positive probability \( p(k) \) which
does not depend on \( j \). So for a random spectrum \( H \) the condition
of the proposition is true almost surely.

\subsection{{}}

Now we pass to the proof of the proposition. Denote \( I(k):=\{sl(k)+q\, :\, 
0<|s|<k,\, |q|<k\} \)
and \( M(k)=kl(k)+k \). Clearly (passing to a subsequence if necessary)
we may suppose that \[
I(k+1)\cap [-M(k),M(k)]=\emptyset \ .\]
Let \( f\in L^{0}(\RR ) \) be given. We shall define by induction
an increasing sequence \( \{k_{j}\} \) and {}``blocks'' of coefficients
\( \{c(n)\} \), \( n\in I(k_{j}) \); all other coefficients of the
expansion (\ref{olv2}) will be zero. We denote:\begin{eqnarray*}
A_{j} & := & \sum _{I(k_{j})}c(n)\exp (i\lambda(n)x)\\
S_{N} & := & \sum _{j\leq N}A_{j}\quad N\in \ZZ ^{+},\, S_{0}:=0.
\end{eqnarray*}
Fix \( N \) and suppose that the polynomials \( A_{j} \) are already
defined for \( j<N \). Let us describe the \( N \)'th step of the
induction. Set\[
R_{N}:=f-S_{N-1}\ .\]
We need the following result proved in \cite{olevskii}:

 if \( 0<|r(q)|=o(1) 
\)
then the system of exponentials \( \exp i(q+r(q)) \) is complete
in \( L^{0}(\RR ) \), that is, the set of linear combinations is dense
with respect to convergence a.e.

 Using this we find a polynomial\[
F_{N}(x)=\sum a_{q}\exp i(q+2^{-|q|+1})x\]
so that\begin{equation}
\label{olv6}
\mathbf{m}\left\{ x\in [-N\pi ,N\pi ]\, :\, 
|F_{N}(x)-R_{N}(x)|>\frac{1}{N^{4}}\right\} <\frac{1}{N^{2}}\ .
\end{equation}
Next we use Lemma \ref{Lemma_Korner} with\begin{equation}
\label{olv7}
\delta =\delta _{N}=\frac{1}{N^{4}||\widehat{F_{N}}||_{1}},\quad \epsilon 
=\frac{1}{N^{3}}
\end{equation}
and find the corresponding polynomial \( Q_{N} \). Fix a number \( k_{N} \)
large enough:\begin{equation}
\label{olv8}
k_{N}>k_{N-1},\, 3\deg F_{N},\, \frac{||\widehat{Q_{N}}||_{1}}{\delta 
_{N}}
\end{equation}
and set \begin{equation}
\label{olv9}
H_{N}:=F_{N}\cdot (Q_{N})_{[l(k_{N})]}\ .
\end{equation}
One can easily see that \[
\spec H\subset \left\{ sl(k_{N})+q+2^{-|q|+1}\, :\, 0<|s|<k,\, |q|<k\right\} 
\ ,\]
so we can write \begin{equation}
\label{olv10}
H_{N}:=\sum _{\substack {0<|s|<k_{N}\brk |q|<k_{N}}}b(N,s,q)\exp 
i(sl(k_{N})+q+2^{-|q|+1})x\ .
\end{equation}
Finally we set:\begin{equation}
\label{olv11}
A_{N}:=\sum _{\substack {0<|s|<k_{N}\brk |q|<k_{N}}}b(N,s,q)\exp 
i\lambda(sl(k_{N})+q)x\equiv \sum
_{I(k_{N})}c(n)\exp i\lambda(n)x\ .
\end{equation}

\subsection{{}}

Now we show that \begin{equation}
\label{olv12}
S_{N}\rightarrow f\quad \textrm{a}.\textrm{e}.
\end{equation}
For this first we get from \ref{req_Q1_lem_korner} of Lemma 
\ref{Lemma_Korner}:
\begin{eqnarray}
\lefteqn {\mathbf{m}\left( \left\{ x\in [-N\pi ,N\pi ]\: :\: \left| 
H_{N}(x)-F_{N}(x)\right| >\frac{1}{N^{4}}\right\} \right) \leq } &  &
\nonumber \\
& \qquad  & \leq N\cdot \mathbf{m}\left( \left\{ x\in [-\pi ,\pi ]\, :\, 
|Q_{N}-1|\geq \delta _{N}\right\} \right) =O(N^{-2})\quad .\label{olv13}
\end{eqnarray}
Further, (\ref{olv9}), (\ref{olv8}) and (\ref{olv7}) imply: \begin{equation}
\label{Hn_hat_1}
||\widehat{H_{N}}||_{1}=||\widehat{F_{N}}||_{1}\cdot 
||\widehat{Q_{N}}||_{1}<||\widehat{F_{N}}||_{1}k_{N}\delta 
_{N}=\frac{k_{N}}{N^{4}}\ ,
\end{equation}
so we can estimate, using (\ref{olv10}), (\ref{olv11}), (\ref{olv5}),
(\ref{Hn_hat_1}), (\ref{olv3}) and (\ref{olv8}):
\begin{eqnarray}
\lefteqn {||A_{N}-H_{N}||_{L^{\infty }(-\pi N,\pi N)}\leq } &  & \nonumber 
\\
& \qquad \leq  & ||\widehat{H_{N}}||_{1}\cdot \max _{\substack 
{0<|s|<k_{N}\brk |q|<k_{N}}}||\exp i(r(sl(k_{N})+q)-2^{-|q|+1})x-1||_{\infty 
}\nonumber \\
& \qquad < & \frac{k_{N}}{N^{4}}\cdot \frac{1}{k^{2}_{N}}\cdot \pi 
N=O(N^{-4})\ .\label{olv14}
\end{eqnarray}
Finally, we have from (\ref{olv6}), (\ref{olv13}) and (\ref{olv14}):\[
\mathbf{m}\left\{ x\in [-N\pi ,N\pi ]\, :\, 
|A_{N}(x)-R_{N}(x)|>\frac{C}{N^{4}}\right\}
=O(N^{-2})\ ,\]
so\begin{equation}
\label{olv15}
R_{N+1}=R_{N}-A_{N}=O(N^{-4})\quad \textrm{a}.\textrm{e}.,
\end{equation}
and (\ref{olv12}) follows.

\subsection{{}}

At last:\begin{equation}
\label{olv16}
A^{*}_{N}\rightarrow 0\quad \textrm{a}.\textrm{e}.
\end{equation}
Indeed, estimating as in (\ref{olv14}) we see: \[
A^{*}_{N}<H^{*}_{N}+O(N^{-4})\quad \textrm{a}.\textrm{e}.\]
Since \( l(k_{N})>k_{N}>3\deg F_{N} \) we can use (\ref{olv4}) and
get:\[
H^{*}_{N}(x)<|F_{N}(x)|\cdot ||Q_{N}^{*}||_{L^{\infty }(-\pi ,\pi 
)}+2||\widehat{Q}||_{\infty }\cdot
||\widehat{F}||_{1}\ .\]
The first term on the right hand side is \( O(N^{-1}) \) a.e. due
to (\ref{olv6}), (\ref{olv15}) and \ref{req_Sq_lem_korner}. The
last term is \( O(N^{-4}) \) because of (i) and (\ref{olv7}). Clearly
(\ref{olv12}) and (\ref{olv16}) imply the decomposition (\ref{olv2})
and this completes the proof.\hfill $\square$

\section{Remarks}

\subsection{{}}

One can see that the result holds for \( r(n) \) uniformly
distributed
on any fixed neighbourhood of zero. Moreover, it holds true for \( r(n) \)
uniformly distributed on \( [-d _{n},d _{n}] \) if \( d _{n} \)
decrease slowly enough. This allows one to cover in full generality the
result from \cite{KO-99} where a Menshov spectrum \( \{n+o(1)\} \)
has been constructed. But \( d _{n} \) really must decay slowly.
In contrast to the completeness property which
occurs for 
\textbf{any}
(nonzero) perturbation \( r(n)=o(1) \) (the result from \cite{olevskii}
used above), the  following simple observation is true:
\theoremstyle{nopoint}
\newtheorem*{ialt}{}
\begin{ialt} If \( \lambda(n)=n+O(n^{-\alpha }) \), \( \alpha >0 \), then
it is not a Menshov spectrum. \end{ialt} Indeed, if \( \Delta 
f:=f(x)-f(x-2\pi ) \)
then \[
\Delta ^{k}\left( \sum c_{n}e^{i\lambda_{n}x}\right) =\sum c_{n}O(n^{-\alpha 
k})e^{i\lambda_{n}x}\ .\]
Hence for \( k \) sufficiently large (depending only on \( \alpha  \))
the representation (2) implies that \( \Delta ^{k}f \) equals  a smooth
function  a.e., so \({ \lambda(n)} \) cannot be a
Menshov spectrum for \( \RR  \).

\subsection{{}}

In \cite{KO01} we studied Menshov spectra in the periodic case. The
main results of that paper can be extended to Menshov spectra in \( \RR  \).
For example, Menshov spectra for \( \RR  \) may be quite sparse,
up to {}``almost hadamarian lacunarity''. More precisely:

 {\it for any
\( \epsilon (n) \) decreasing to zero one can construct a (symmetric)
Menshov spectrum \( H \) for \( \RR  \) such that
 \(\lambda(n+1)/\lambda(n)>1+\epsilon (n)\),  \( n\in \ZZ ^{+}
\).}

 This is an analog of Theorem 1 from \cite{KO01}
and the proof is basically the same.

\bigskip
\medskip

\noindent
School of Mathematical Sciences, Tel Aviv University, 
Tel Aviv 69978, Israel\\
E-mail addresses: gady@post.tau.ac.il,
olevskii@post.tau.ac.il

\end{document}